\documentclass{article} 

\usepackage{amsmath,amsthm}     
\usepackage{graphicx}     
\usepackage{hyperref} 
\usepackage{url}
\usepackage{amsfonts} 

\usepackage{multicol}

\theoremstyle{theorem}

\theoremstyle{definition}

\allowdisplaybreaks

\makeatletter
\@addtoreset{footnote}{page}
\makeatother

\begin{document}

\title{Elimination and Factorization}

\author{Gilbert Strang\\               
Massachusetts Institute of Technology\\    
gilstrang@gmail.com
\date{}\\
}                      

\maketitle

\vspace{-11pt}
\begin{center}
MSC Classifications\quad 15A21\quad   15A23\quad 97\qquad\;\; ORCID\quad 0000-0001-8511-7862
\end{center}

\vspace{12pt}
\begin{tabular}{@{\hspace{-2pt}}c}
\textbf{{\fontsize{10}{42}\selectfont \textbf{Abstract}}}\\
\parbox{11.2cm}{If a matrix $A$ has rank $r$, then its row echelon form (from elimination) contains the identity matrix in its first $r$ independent columns. How do we \emph{interpret the matrix} $F$ that appears in the remaining columns of that echelon form\,? $F$ multiplies those first $r$ independent columns of $A$ to give its $n-r$ dependent columns. Then $F$ reveals bases for the row space and the nullspace of the original matrix $A$. And $F$ is the key to the column-row factorization $\boldsymbol{A}=\boldsymbol{CR}$.}
\end{tabular}

\vspace{20pt}
\noindent
$\boldsymbol{1.}$ Elimination must be just about the oldest algorithm in linear algebra. By systematically producing zeros in a matrix, it simplifies the solution of $m$ equations $A\boldsymbol{x}=\boldsymbol{b}$. We take as example this $3$ by $4$ matrix $A$, with row $1\,+$ row $2=$ row $3$. Then its rank is $r=2$, and we execute the familiar elimination steps to find its \emph{reduced row echelon form} $Z$\,:

\vspace{3pt}
$$
\boldsymbol{A}=\left[\begin{array}{rrrr}1&2&11&17\\3&7&37&57\\4&9&48&74\\\end{array}\right]\rightarrow \boldsymbol{Z}=\left[\begin{array}{rrrr}1&0&3&5\\0&1&4&6\\0&0&0&0\\\end{array}\right].
$$

\vspace{7pt}
At this point, we pause the algorithm to ask a question\,: \textbf{How is $\boldsymbol{Z}$ related to $\boldsymbol{A}$\,?} One answer comes from the fundamental subspaces associated with $\boldsymbol{A}$\,:

\begin{enumerate}
	\item [1)] The two nonzero rows of $\boldsymbol{Z}$ (call them $\boldsymbol{R}$) are a basis for the row space of $\boldsymbol{A}$.
	\item [2)] The first two columns of $\boldsymbol{A}$ (call them $\boldsymbol{C}$) are a basis for the column space of $\boldsymbol{A}$.
	\item [3)] The nullspace of $\boldsymbol{Z}$ equals the nullspace of $\boldsymbol{A}$ (orthogonal to the same row space).
\end{enumerate}

\noindent
Those were our reasons for elimination in the first place. ``Simplify the matrix $A$ without losing the information it contains.'' By applying the same steps to the right hand side of $A\boldsymbol{x}=\boldsymbol{b}$, we reach an equation $Z\boldsymbol{x}=\boldsymbol{d}$---with the same solutions $\boldsymbol{x}$ and the simpler matrix $Z$.

The object of this short note is to express the result of elimination in a different way. This factorization cannot be new, but it deserves new emphasis.

\vspace{10pt}
\noindent
\textbf{Elimination factors $\boldsymbol{A}$ into $\boldsymbol{C}$ times $\boldsymbol{R=(m\times r)}$ times $\boldsymbol{(r\times n)}$}

\vspace{4pt}
$$
\boldsymbol{A}=\left[\begin{array}{rrrr}1&2&11&17\\3&7&37&57\\4&9&48&74\\\end{array}\right]=\left[\begin{array}{rr}1&2\\3&7\\4&9\\\end{array}\right]\begin{array}{@{\!}r@{\!}}\left[\begin{array}{rrrr}1&0&3&5\\0&1&4&6\\\end{array}\right]\\\\\end{array}=\boldsymbol{CR}
$$

\vspace{6pt}
\noindent
$C$ has full column rank $r=2$ and $R$ has full row rank $r=2$. When we establish that $A=CR$ is true for every matrix $A$, this factorization brings with it a proof of the first great theorem in linear algebra\,: \textbf{Column rank equals row rank}.

\vspace{11pt}
\noindent
$\boldsymbol{2.}$\quad Here is a description of $C$ and $R$ that is independent of the algorithm (row operations) that computes them.

\vspace{5pt}
\textbf{Suppose the first $\boldsymbol{r}$ independent columns of $\boldsymbol{A}$ go into $\boldsymbol{C}$}. \textbf{Then the other $\boldsymbol{n-r}$ columns of $\boldsymbol{A}$ must be combinations $\boldsymbol{CF}$ of those independent columns in $\boldsymbol{C}$. That key matrix $\boldsymbol{F}$ is part of the row factor $\boldsymbol{R=\left[\begin{array}{@{\;}rr@{\;}}I&F\end{array}\right]P}$, with $\boldsymbol{r}$ independent rows.} Then right away we have $A=CR$\,:

\vspace{-11pt}
\begin{equation}
\!\boldsymbol{A}\!=\!\boldsymbol{CR}=\!\left[\begin{array}{@{\:}rr@{\:}}\boldsymbol{C}&\boldsymbol{CF}\end{array}\right]\boldsymbol{P}\!=\!\left[\begin{array}{@{\:}rr@{\:}}\textbf{Indep cols}&\textbf{Dep cols}\end{array}\right]\textbf{Permute cols}\!
\label{eq:appelim01}
\end{equation}

\vspace{11pt}
\noindent
If the $r$ independent columns come first in $A$, that permutation matrix will be $P=I$. Otherwise we need $P$ to permute the $n$ columns of $C$ and $CF$ into correct position in $A$. Here is an example of $A=\!C\!\left[\begin{array}{@{\,}rr@{\,}}I&F\\\end{array}\right]\!P$ in which $P$ exchanges columns $2$ and $3$\,:

\vspace{-8pt}
\begin{equation}
\!\!\boldsymbol{A}\!=\!\!\left[\begin{array}{@{\,}rrrr@{\,}}1&2&3&4\\1&2&4&5\\\end{array}\right]\!\!=\!\!\left[\begin{array}{@{\,}rr@{\,}}1&3\\1&4\\\end{array}\right]\!\!\left[\begin{array}{@{\,}rrrr@{\,}}1&2&0&1\\0&0&1&1\\\end{array}\right]\!\!=\!\!\left[\begin{array}{@{\,}rr@{\,}}1&3\\1&4\\\end{array}\right]\!\!\left[\begin{array}{@{\,}rrrr@{\,}}1&0&2&1\\0&1&0&1\\\end{array}\right]\!P\!=\!\boldsymbol{CR}.
\label{eq:appelim02}
\end{equation}

\vspace{-8pt}
\vspace{11pt}
\noindent
The essential information in $\hbox{\textbf{rref}}(A)$ is the list of $r$ independent columns of $A$, and the matrix $F$ ($r$ by $n-r$) that combines those independent columns to give the $n-r$ dependent columns $CF$ in $A$. This uniquely defines $\hbox{\textbf{rref}}(A)$.

\vspace{11pt}
\noindent
$\boldsymbol{3.}$\quad The factorization $A=CR$ is confirmed. But how do we determine the first $r$ independent columns in $A$ (going into $C$) and the dependencies of the remaining $n-r$ columns $CF$\,? This is the moment for \textbf{row operations} on $A$. Three operations are allowed, to put $A$ into its reduced row echelon form $Z=\textbf{rref}(A)$\,:

\begin{enumerate}
	\item[(a)] Subtract a multiple of one row from another row (below or above)
	\item[(b)] Exchange two rows
	\item[(c)] Divide a row by its first nonzero entry
\end{enumerate}

\newpage
\noindent
All linear algebra teachers and a positive fraction of students know those steps and their outcome $\textbf{rref}(A)$. It contains an $r$ by $r$ identity matrix $I$ (only zeros can precede those $1$'s). The position of $I$ reveals the first $r$ independent columns of $A$. And equation (\ref{eq:appelim01}) above reveals the meaning of $F$\,! It tells us how the $n-r$ \emph{dependent} columns $CF$ of $A$ come from the independent columns in $C$. The remaining $m-r$ dependent rows of $A$ must become zero rows in $Z$\,:

\vspace{-5pt}
\begin{equation}
\textbf{Elimination reduces\;}  \boldsymbol{A} \textbf{ \;to\; } \boldsymbol{Z}=\hbox{\textbf{rref}}(A)=\left[\begin{array}{cc}\boldsymbol{I}&\boldsymbol{F}\\\boldsymbol{0}&\boldsymbol{0}\\\end{array}\right]\boldsymbol{P}
\label{eq:appelim03}
\end{equation}

\noindent
All our row operations (a)(b)(c) are invertible. (This is Gauss-Jordan elimination\,: operating on rows above the pivot row as well as below.) But the matrix that reduces $A$ to this echelon form is less important than the \textbf{factorization} $\boldsymbol{A=CR}$ that it uncovers in equation (\ref{eq:appelim01}).

\vspace{11pt}
\noindent
$\boldsymbol{4.}$\quad Before we apply $A=CR$ to solving $A\boldsymbol{x}=\boldsymbol{0}$, we want to give a column by column (left to right) construction of $\textbf{rref}(A)$ from $A$. After elimination on $k$ columns, that part of the matrix is in its own $\textbf{rref}$ form. We are ready for elimination on the current column $k+1$. This new column has an upper part $\boldsymbol{u}$ and a lower part $\boldsymbol{\ell}$\,:

\vspace{-14pt}
\begin{equation}
\textbf{First }\boldsymbol{k+1} \textbf{ columns }\qquad \left[\begin{array}{cc}I_k&F_k\\0&0\\\end{array}\right]P_k\text{ \;followed by\; } \left[\begin{array}{c}\boldsymbol{u}\\\boldsymbol{\ell}\\\end{array}\right].
\label{eq:}
\end{equation}

\vspace{1pt}
The big question is\,: \textbf{Does this new column $\boldsymbol{k+1}$ join with $\boldsymbol{I_k}$ or $\boldsymbol{F_k}$\,?}

\vspace{5pt}
\textbf{If $\boldsymbol{\ell}$ is all zeros}, the new column is \textbf{dependent} on the first $k$ columns. Then $\boldsymbol{u}$ joins with $F_k$ to produce $F_{k+1}$ in the next step to column $k+2$.

\vspace{5pt}
\textbf{If $\boldsymbol{\ell}$ is not all zero}, the new column is \textbf{independent} of the first $k$ columns. Pick any nonzero in $\boldsymbol{\ell}$ (preferably the largest) as the pivot. Move that row of $A$ to the top of $\boldsymbol{\ell}$. Then use that pivot row to zero out (by standard elimination) all the rest of column $k+1$. (That step is expected to change the columns after $k+1$.) Column $k+1$ joins with $I_k$ to produce $I_{k+1}$. We adjust $P_k$ and we  are ready for column $k+2$.

\vspace{5pt}
At the end of elimination, we have a most desirable list of column numbers. They tell us the \textbf{first $\boldsymbol{r}$ independent columns of $\boldsymbol{A}$}. Those are the columns of $C$. They led to the identity matrix $I_{r\text{ by } r}$ in the row factor $R$ of $A=CR$.

\vspace{11pt}
\noindent
$\boldsymbol{5.}$\quad What is achieved by reducing $A$ to $\textbf{rref}(A)$\,? The row space is not changed\,! Then its orthogonal complement (\textbf{the nullspace of} $\boldsymbol{A}$) is not changed. Each column of $CF$ tells us how a dependent column of $A$ is a combinaton of the independent columns in $C$. Effectively, \textbf{the columns of $\boldsymbol{F}$ are telling us $\boldsymbol{n-r}$ solutions to $\boldsymbol{Ax=0}$}.

This is easiest to see by the example in Section $1$\,: 

\begin{center}
$\begin{array}{r@{\,}r@{\,}r@{\,}r@{\,}r@{\,}r@{\,}r@{\,}r}x_1&+&2x_2&+&11x_3&+&17x_4&=0\\3x_1&+&7x_2&+&37x_3&+&57x_4&=0\\\end{array}$\quad reduces to\quad 
$\begin{array}{r@{\,}r@{\,}r@{\,}r@{\,}r@{\,}r}x_1&+&3x_3&+&5x_4&=0\\x_2&+&4x_3&+&6x_4&=0\\\end{array}$
\end{center}

\vspace{5pt}
\noindent
The solution with $x_3=1$ and $x_4=0$ is $\boldsymbol{x}=\left[\begin{array}{@{\,}rrrr}-3&-4&1&0\\\end{array}\right]^{\mathrm{T}}$. Notice $3$ and $4$ from $F$. The second solution with $x_3=0$ and $x_4=1$ is $\boldsymbol{x}=\left[\begin{array}{@{\,}rrrr}-5&-6&0&1\\\end{array}\right]^{\mathrm{T}}$. Those solutions are the two columns of $\boldsymbol{X}$ in $\boldsymbol{AX=0}$. (This example has $P=I$.) \textbf{The $\boldsymbol{n-r}$ columns of $\boldsymbol{X}$ are a natural basis for the nullspace of} $\boldsymbol{A}$\,:

\vspace{-11pt}

$$\boldsymbol{A}\hspace{-1pt}=\hspace{-1pt}\boldsymbol{C}\left[\begin{array}{@{\;}rr@{\;}}\boldsymbol{I}&\boldsymbol{F}\\\end{array}\right]\hspace{-1pt}\boldsymbol{P} \text{ \;multiplies\; } \boldsymbol{X}\hspace{-1pt}=\hspace{-1pt}\boldsymbol{P^{\mathrm{T}}}\left[\begin{array}{@{\;}c@{\;}}\boldsymbol{-F}\\\boldsymbol{I_{n-r}}\\\end{array}\right] \text{ to give\, }$$

\begin{equation}
\hspace{1pt}\boldsymbol{AX=-CF+CF=0}.
\label{eq:}
\end{equation}

\vspace{11pt}
\noindent
With $PP^{\mathrm{T}}=I$ for the permutation, each column of $X$ solves $A\boldsymbol{x}=\boldsymbol{0}$. Those\linebreak $n-r$ solutions in $X$ tell us what we know\,: Each dependent column of $A$ is a combination of the independent columns in $C$. In the example, column $3=3\,$(column $1$) $+\;4\,$(column $2$).

\vspace{5pt}
Gauss-Jordan elimination leading to $A=CR$ is less efficient than the Gauss process that directly solves $A\boldsymbol{x}=\boldsymbol{b}$. The latter stops at a triangular system $U\boldsymbol{x}=\boldsymbol{c}$\,: back substitution produces $\boldsymbol{x}$. Gauss-Jordan has the extra cost of eliminating upwards. If we only want to solve equations, stopping at a triangular factorization is faster.

\vspace{11pt}
\noindent
$\boldsymbol{6.}$\quad \textbf{Block elimination}\quad The row operations that reduce $A$ to its echelon form produce one zero at a time. A key part of that echelon form is an $r$ by $r$ identity matrix. If we think on a larger scale---instead of one row at a time---\textbf{that output $\boldsymbol{I}$ tells us that some $\boldsymbol{r}$ by $\boldsymbol{r}$ matrix has been inverted}. Following that lead brings a ``matrix understanding'' of elimination.

\vspace{5pt}
Suppose that the matrix $\boldsymbol{W}$ in the first $r$ rows and columns of $A$ is invertible. Then \emph{elimination takes all its instructions from} $W$\,! One entry at a time---or all at once by ``block elimination''---\textbf{$\boldsymbol{W}$ will change to $\boldsymbol{I}$}. In other words, the first $r$ rows of $A$ will yield $I$ and $F$. \textbf{This identifies $\boldsymbol{F}$ as $\boldsymbol{W^{-1}H}$}. And the last $m-r$ rows will become zero rows.

\vspace{-9pt}
\begin{equation}
\!\textbf{Block elimination}\quad \boldsymbol{A}=\!\left[\begin{array}{cc}\boldsymbol{W}&\boldsymbol{H}\\\boldsymbol{J}&\boldsymbol{K}\\\end{array}\right]\text{ reduces to } \left[\begin{array}{cc}\boldsymbol{I}&\boldsymbol{F}\\\boldsymbol{0}&\boldsymbol{0}\\\end{array}\right]\!=\textbf{rref}\boldsymbol{(A)}.
\label{eq:}
\end{equation}

\vspace{4pt}
\noindent
That just expresses the facts of linear algebra\,: If $A$ begins with $r$ independent rows and its rank is $r$, then the remaining $m-r$ rows are combinations of those first $r$ rows \,:\! $\left[\begin{array}{@{\:}rr@{\:}}J&K\\\end{array}\right]\!=\!JW^{-1}\left[\begin{array}{@{\:}rr@{\:}}W&H\\\end{array}\right]$. Those $m-r$ rows become zero rows in elimination.

\vspace{5pt}
In general that first $r$ by $r$ block might not be invertible. But elimination will find $W$. We can move $W$ to the upper left corner by row and column permutations $\boldsymbol{P_{\text{\small{$r$}}}}$ and $\boldsymbol{P_{\text{\small{$c$}}}}$.

Then the full expression of block elimination to reduced row echelon form is

\vspace{-1pt}
\begin{equation}
\boldsymbol{P_{\text{\small{$r$}}}AP_{\text{\small{$c$}}}=\left[\begin{array}{cc}W&H\\J&K\\\end{array}\right]\rightarrow\left[\begin{array}{cc}I&W^{-1}H\\0&0\\\end{array}\right]}
\label{eq:}
\end{equation}

\vspace{9pt}
\noindent
$\boldsymbol{7.}$\quad This raises an interesting point. Since $A$ has rank $r$, we know that it has $r$ independent rows and $r$ independent columns. Suppose those rows are in a submatrix $B$ and those columns are in a submatrix $C$. Is it always true that the $r$ by $r$ ``intersection'' $W$ of those rows with those columns will be \textbf{invertible}\,? (Everyone agrees that somewhere in $A$ there is an $r$ by $r$ invertible submatrix. The question is whether $B\cap C$ can be counted on to provide such a submatrix.) Is $W$ automatically full rank\,?

The answer is \emph{yes}. \textbf{The intersection of $\boldsymbol{r}$ independent rows of $\boldsymbol{A}$ with $\boldsymbol{r}$ independent columns does produce a matrix $\boldsymbol{W}$ of rank $\boldsymbol{r}$. $\boldsymbol{W}$ is invertible.}

\vspace{10pt}
\noindent
\textbf{Proof} \cite{3}\,: Every column of $A$ is a combination of the $r$ columns of $C$.

\noindent
Every column of the $r$ by $n$ matrix $B$ is the same combination of the $r$ columns of $W$.

\noindent
Since $B$ has rank $r$, its column space is all of $\hbox{\textbf{R}}^{\boldsymbol{r}}$.

\noindent
Then the column space of $W$ is also $\hbox{\textbf{R}}^{\boldsymbol{r}}$  and \textbf{the square submatrix $\boldsymbol{W}$ has rank $\boldsymbol{r}$}.

\vspace{10pt}
\noindent
\mbox{\textbf{Key words}\,: Column-row factorization, $A=CR$, Elimination, Echelon form $R$}

\end{document}